\documentclass[12pt,oneside,english]{amsart}
\usepackage[latin1]{inputenc}
\usepackage{amssymb}

\makeatletter
\newtheorem{theorem}{Theorem}
\newtheorem{lemma}{Lemma}

\newtheorem{definition}{Definition}

\newtheorem{remark}{Remark}

\newtheorem{variation}{Variation}
\newtheorem{tab}{Table}
\numberwithin{equation}{section}
\usepackage{babel}

\makeatother
\begin{document}
\baselineskip=17pt

\title[Enumeration of two-color bracelets with variations]{A problem of enumeration of two-color bracelets with several variations}

\author{Vladimir Shevelev}
\address{Department of Mathematics \\Ben-Gurion University of the
 Negev\\Beer-Sheva 84105, Israel. e-mail:shevelev@bgu.ac.il}
\subjclass{05A15}

\begin{abstract}
We consider the problem of enumeration of incongruent two-color bracelets of $n$ beads, $k$ of which are black,
and study several natural variations of this problem. We also give recursion formulas for enumeration of $t$-color bracelets, $t\geq3.$

\end{abstract}

\maketitle

\section{Introduction}
     Professor Richard H.Reis (South-East University of Massachusetts,
USA) in 1978 put the problem: "Let a circumference is split by the
same $n$ parts. It is required to find the number $R(n,k)$ of the
incongruent convex k-gons,which could be obtained by connection of
some $k$ from $n$ dividing points. Two k-gons are considered
congruent if they are coincided at the rotation of one relatively
other along the circumference and (or) by reflection of one of the
k-gons relatively some diameter."

     In 1979 Hansraj Gupta \cite{1}gave the solution of the Reis
problem.

\begin{theorem}\label{t1}(H.Gupta)
\begin{equation}\label{1.1}
R(n,k)=\frac 1
2\left(\binom{\lfloor\frac{n-h_k}{2}\rfloor}{\lfloor\frac{k}{2}\rfloor}+\frac
1 k\sum_{d|(k,n)}\varphi(d)\binom{\frac n d -1}{\frac k d
-1}\right),
\end{equation}
where $h_k\equiv k(\mod{2}), \;\; h_k=0$  or 1,\;\;$(n,k)$ is $\gcd (n,k),\enskip \varphi(n)$-the
Euler function.
\end{theorem}

     Consider some convex polygon with the tops in the circumference
splitting points, "1" or "0" is put in accordance to each splitting
point depending on whether a top of the polygon is in the point.
Thus, there is the mutual one-to-one correspondence between the  set
of convex polygons with the tops in the circumference splitting
points and the set of all $(0,1)$-configurations with the elements
in these points.

     Using this bijection and calculating the cyclic index of the
Dihedral group appearing here, the author \cite{4} gave a short proof
of Theorem 1. Besides, this bijection shows that formula (\ref{1.1}) solves the
problem of enumeration of two-color bracelets of $n$ beads, $k$ of which are black and $n-k$ are white.
In turn, it allows to obtain several simple formulas for series of sequences in \cite{6}
(see, e.g., the author's explicit formulas for sequences A032279-A032282,\enskip A005513-A005516). \newline
\indent  Note also that, quite recently, the author found an application of enumeration of two-color bracelets
to some questions of the permanent theory (see \cite{5}, Section 5).\newline

\indent Let $n$ beads of a bracelets are located in $n$ dividing points of a circumference which is split by the
same $n$ parts. Let $\mathcal{T}$ be cyclic group of turns with the generating element $\tau=e^{\frac{2\pi i}{n}}.$

\begin{definition}\label{d1}
We call two-color bracelet of $n$ beads  \upshape symmetric respectively rotation (two-color SR-bracelet)\slshape \enskip if its cyclic grope of turns is a proper subgroup of $\mathcal{T}.$
\end{definition}
\begin{remark}
When we say about a two-color bracelet, we exclude cases when it contains only white (black) color.
\end{remark}
\indent Consider several variations of enumeration of two-color bracelets.

\begin{variation} \label{V1}
To find the number $N_{n}$ of all incongruent two-color SR-bracelets.
\end{variation}
\begin{variation} \label{V2} To find the
number $N^{(k)}(n)$ of those two-color SR-bracelets
which have exactly $k$ black beads.
\end{variation}
\begin{variation} \label{V3} To find the
number $S_{n}$ of those two-color SR-bracelets
which have a diameter of symmetry.
\end{variation}
\begin{variation} \label{V4} To find the number $S^{(k)}(n)$ of those two-color SR-bracelets
which have exactly $k$ black beads and a diameter of symmetry.
\end{variation}
\begin{variation} \label{V5} Let $m$ be a positive
integer. To find the number $N_{n,m}$ of all incongruent two-color SR-bracelets
with isolated black beads such that between every two black beads there exist at least $m$ white ones.
\end{variation}

\begin{variation} \label{V6} For $m\geq 1$, to find the
number $S_{n,m}$ of those two-color SR-bracelets
which have a diameter of symmetry.
\end{variation}

\begin{variation} \label{V7} For $m\geq 1,$ to find the
number $N^{(k)}_{n,m}$ of those two-color SR-bracelets in Variation $\ref{V2}$
which have exactly $k$ black beads.
\end{variation}

\begin{variation} \label{V8} For $m\geq 1,$ to find the
number $S^{(k)}_{n,m}$ of those two-color SR-bracelets in Variation $\ref{V3}$
which have exactly $k$ black beads.
\end{variation}
Notice that, $N_{n,m}-S_{n,m}$ is the number of those two-color SR-bracelets,\newpage
  none of which has a diameter of symmetry; $N_{n,m}^{(k)}-S_{n,m}^{(k)}$ is
the same two-color SR-bracelets having exactly $k$ 1`s.\newline
\indent Some words about structure of the article. Section 2 is
devoted to solutions of Variations \ref{V1}-\ref{V4}. In Section 3 we introduce two different
generalizations of the Fibonacci numbers. In Section 4 we solve
 Variations \ref{V5}-\ref{V8}. In Section 5 we consider an interesting example and the connected
with it some numerical results. In conclusion, in  Section 6 we discuss
two open questions (Variation \ref{V9}-\ref{V10}) and give enumeration of $t$-color SR-bracelets (Variation \ref{V11}-\ref{V12}) and Theorems \ref{t4}-\ref{t6}.

\section{Variations 1-4}

\begin{theorem}\label{t2}
The following formulas take place:

\begin{equation}\label{2.1}
N_n=-2-\sum\limits_{d|n,\enskip d\geq 2}\mu(d)\alpha_{\frac n d},
 \end{equation}
\begin{equation}\label{2.2}
S_n=-2-\sum\limits_{d|n,\enskip d\geq 2}\mu(d)\beta_{\frac n d},
\end{equation}
\begin{equation}\label{2.3}
N^{(k)}_n=-\sum\limits_{d|(k,n)\enskip d\geq 2}\mu(d)R\left(\frac n
d,\frac k d\right),
\end{equation}
\begin{equation}\label{2.4}
S^{(k)}_n=-\sum\limits_{d|(k,n),\enskip d\geq 2}\mu(d)R^1\left(\frac n
d,\frac k d\right),
\end{equation}
where $\mu(n)$ is the Mobius function,
\begin{equation}\label{2.5}
\beta_n=(5+(-1)^n)2^{\lfloor\frac{n-3}{2}\rfloor},
\end{equation}

\begin{equation}\label{2.6}
\alpha_n=\frac 1 n \sum\limits_{d|n}\varphi(d) 2^{\frac n d -1}+\frac{\beta_n}{2},
\end{equation}

\begin{equation}\label{2.7}
R^1(n,k)=\begin{pmatrix}\left\lfloor\frac{n-h_k}{2}\right\rfloor\\
\left\lfloor\frac k 2\right\rfloor\end{pmatrix}
\end{equation}
and $R(n,k)$  is defined by $(\ref{1.1}).$
\end{theorem}

\bfseries Proof.\mdseries\;\;\; 1) Summing (\ref{1.1}) by $k$ from 1 to
$n-1$ we find the number $\lambda_n$ of all incongruent
of all incongruent two-color bracelets:
\begin{equation}\label{2.8}
\lambda_n=\sum\limits^{n-1}_{k=1} R(n,k).
\end{equation}

Let $d, 1\leq d \leq n,$ be a divisor of $n$. Denote via $\nu_d=\nu_d(n)$ the
number of\newpage incongruent bracelets with the minimal angle of
self-coincidence equals to $\frac{2\pi}{n}d$.
Then
\begin{equation}\label{2.9}
\lambda_n=\sum\limits_{d|n}\nu_d
\end{equation}

and, by definition of $N_n,$ we have

\begin{equation}\label{2.10}
N_n= \lambda_n-\nu_1-\nu_n.
\end{equation}

Using the Mobius inverse formula (cf., e.g., \cite{3}), we find from
(\ref{2.9})
$$
\nu_n=\sum\limits_{d|n}\mu(d) \lambda_{\frac n
d}=\lambda_n+\sum\limits_{d|n,\enskip d\geq 2}\mu(d) \lambda_{\frac n d}.
$$

Now (\ref{2.8}) implies
\begin{equation}\label{2.11}
N_n=-\sum\limits_{d|n,\enskip 2\leq d\leq \frac{n}{2}}\mu(d) \lambda_{\frac n d}.
\end{equation}

To complete the proof of (\ref{2.1}), we need two
technical lemmas.

\begin{lemma}\label{L1}   For $h_k\equiv k (\mod{2}),\;\;h_k=0 \; or \; 1,$  we have
\begin{equation}\label{2.12}
\sum\limits^{n-1}_{k=1}\begin{pmatrix}\left\lfloor\frac{n-h_k}{2}\right\rfloor\\
\lfloor\frac k
2\rfloor\end{pmatrix}=(5+(-1)^n)2^{\lfloor\frac{n-3}{2}\rfloor}-2.
\end{equation}

\end{lemma}

\bfseries Proof.\mdseries\;\;\;Indeed, for even $n,$ we have
\begin{equation}\label{2.13}
\sum\limits^{n-1}_{k=1}\begin{pmatrix}\left\lfloor\frac{n-h_k}{2}\right\rfloor\\
\lfloor\frac k
2\rfloor\end{pmatrix}=\sum\limits_{k=1,3,\ldots,n-1}\begin{pmatrix}\frac{n-2}{2}\\
\frac{k-1}{2}\end{pmatrix}+\sum\limits_{k=2,4,\ldots,n-2}\begin{pmatrix}\frac {n} 2\\
\frac k 2\end{pmatrix}=
\end{equation}

\begin{equation}\label{2.14}
=\sum\limits^{\frac{n-2}{2}}_{t=0}\begin{pmatrix}\frac{n-2}{2}\\
t\end{pmatrix}+\sum\limits^{\frac {n-2} 2}_{t=1}\begin{pmatrix}\frac n 2\\
t\end{pmatrix}=2^{\frac{n-2}{2}}+2^{\frac n 2}-2= 3\cdot
2^{\frac{n-2}{2}}-2.
\end{equation}

For odd $n,$ we have
\begin{equation}\label{2.15}
\sum\limits^{n-1}_{k=1}\begin{pmatrix}\left\lfloor\frac{n-h_k}{2}\right\rfloor\\
\lfloor\frac k 2\rfloor\end{pmatrix}=\sum\limits_{k=1,3,\ldots,n-2}\begin{pmatrix}
\frac{n-1}{2}\\
\frac{k-1}{2}\end{pmatrix}+\sum\limits_{k=2,4,\ldots,n-1}\begin{pmatrix}\frac{n-1}{2}\\
\frac k 2\end{pmatrix}=
\end{equation}

\begin{equation}\label{2.16}
=\sum\limits^{\frac{n-3}{2}}_{t=0}\begin{pmatrix}\frac{n-1}{2}\\
t\end{pmatrix}+\sum\limits_{t=1}^{\frac{n-1}{2}}\begin{pmatrix}
\frac{n-1}{2}\\t\end{pmatrix}=2\cdot2^{\frac{n-1}{2}}-2,
\end{equation}
and the lemma follows from (\ref{2.14}),(\ref{2.16}).\newline
$\blacksquare$
\begin{lemma}\label{L2}
\begin{equation}\label{2.17}
\sum\limits^{n-1}_{k=1}\frac 1 k
\sum\limits_{d|(n,k)}\varphi(d)\begin{pmatrix}\frac n d -1\\
\frac k d -1\end{pmatrix}=\frac 1 n\sum\limits_{d|n}\varphi(d)
2^{\frac n d}-2.
\end{equation}
\end{lemma}
\newpage
\bfseries Proof.\mdseries\;\;\;First of all, notice that
$$
\frac 1 k\begin{pmatrix}\frac n d -1\\\frac k d
-1\end{pmatrix}=\frac 1 n \cdot\frac{\frac n d}{\frac k d}
\begin{pmatrix}\frac n d -1\\\frac k d
-1\end{pmatrix}=\frac 1 n\begin{pmatrix}\frac n d\\\frac k d
\end{pmatrix}.
$$
Therefore, taking into account that, for $n\geq1,$
$\sum_{d|n}\varphi(d)=n,$
in order to prove (\ref{2.17}), it is sufficient to show that
\begin{equation}\label{2.18}
\sum\limits^{n-1}_{k=1}\sum\limits_{d|(n,k)}\varphi(d)\begin{pmatrix}\frac
n d\\\frac k d\end{pmatrix}=\sum\limits_{d|n}\varphi(d) (2^{\frac n
d}-2).
\end{equation}
Putting $k=dd_1$, we have
$$
\sum\limits^{n-1}_{k=1}\sum\limits_{d|(n,k)}\varphi(d)\begin{pmatrix}\frac
n d\\\frac k d\end{pmatrix}=\sum\limits_{\substack{dd_1\leq
n-1\\d|n}}\varphi(d)\begin{pmatrix}\frac n d\\d_1\end{pmatrix}=
$$

$$
\sum\limits_{\substack{1\leq d\leq
\frac{n}{2}\\d|n}}\varphi(d)\sum\limits^{\frac {n} {d}-1}_{d_1=1}
\begin{pmatrix}\frac n
d\\d_1\end{pmatrix}=\sum\limits_{\substack{1\leq d\leq
\frac{n}{2}\\d|n}}\varphi(d)(2^{\frac n d}-2)=
$$
$$\sum\limits_{d|n}\varphi(d)(2^{\frac n d}-2).$$
This proves (\ref{2.17}). \newline
$\blacksquare$\newline

Now from (\ref{2.8}), (\ref{1.1}) and Lemmas 1, 2 we have
$$
\lambda_n=(5+(-1)^n) 2^{\lfloor\frac
{n-5}{2}\rfloor}-1+$$ $$ \frac 1 n
\sum\limits_{d|n}\varphi(d) (2^{\frac n d -1}-1)=\alpha_n-2.
$$

Finally, from (\ref{2.11}), taking into account, that for $n>1,\;\;\sum\limits_{d|n,\enskip d\geq 2}\mu(d)=-1,$ we find
$$N_n=-\sum\limits_{d|n,\enskip d\geq 2}\mu(d) \lambda_{\frac n d}=$$ $$-\sum\limits_{d|n,\enskip d\geq 2}\mu(d) (\alpha_{\frac n d}-2)=-2-\sum\limits_{d|n,\enskip d\geq 2}\mu(d)\alpha_{\frac n d}.$$
$\blacksquare$ \newline
\indent Note that, as it was expected, if $n$ is prime, then $N_n=0.$\newline\newline
2) As it showed in \cite{1}, $R^1(n,k)$ (\ref{2.7}) gives
the number of those $k$-gons or, by the bijection, those two-color bracelets
having exactly $k$ black beads, that are
symmetric respectively a diameter. Therefore, by the same arguments for $R^1(n,k),$
we obtain (\ref{2.2}).\newline
$\blacksquare$ \newline

3) Let $(k,n)>1$.
For fixed $n,\;k,$ let us consider the function $R(x,y)$ (1) on the set

$$
\left\{\left(\frac {n}{(n,k)}m,\frac{k}{(n,k)}
m\right)/m|(n,k)\right\}.
$$

The constriction of the $R(x,y)$ on this set is a function of $m$.
Denote it by  $R_{n,k}(m)$. Put

\begin{equation}\label{2.19}
R_{n,k}(m)=\sum\limits_{d|m}\nu_d,
\end{equation}

where $\nu_d=\nu_d(n,k)$. In particular, $R_{n,k}((n,k))=R(n,k)$.
Using the Mobius inverse formula, according to (\ref{2.18}) we have

$$
\nu_m=\sum\limits_{d|m}\mu(d) R_{n,k}(\frac m
d)=\sum\limits_{d|m}\mu(d) R\left(\frac{n}{(n,k)}\;\frac m
d,\;\frac{k}{(n,k)}\;\frac m d\right).
$$

In particular, for $m=(n,k)$

\begin{equation}\label{2.20}
\nu_{(n,k)}=\sum\limits_{d|(n,k)}\mu(d) R\left(\frac n d,\;\frac k
d\right),
\end{equation}

and (\ref{2.19}), for $m=(n,k),$ has the form
$$
R(n,k)=\sum\limits_{d|(n,k)}\nu_d.
$$

By the definition, we have now
\begin{equation}\label{2.21}
N^{(k)}_n=R(n,k)-\nu_{(n,k)}
\end{equation}

Finally, from (\ref{2.20})-(\ref{2.21}) we deduce (\ref{2.3})\newline
 $\blacksquare$\newline
 4) From the same arguments it follows that, by replacing in (\ref{2.3})
$R$ by $R^1,$ we obtain $S^{(k)}(n).$ \newline
 $\blacksquare$

    For considerations of the further variations, we need two different generalizations
of Fibonacci numbers.

\section{Two generalizations of Fibonacci numbers}\label{s3}

\bfseries Definition 2.\enskip \mdseries  Let $m$ be a positive integer. We call $m$-Fibonacci numbers of type 1 the
sequence which is defined by the recursion
\begin{equation}\label{3.1}
F^{(m)}_n= F^{(m)}_{n-1} + F^{(m)}_{n-m-1},
\end{equation}
with the initial conditions
\begin{equation}\label{3.2}
F^{(m)}_0= F^{(m)}_1 =\ldots= F^{(m)}_m = 1.
\end{equation}
\newpage
For $m=1$ we obtain the very Fibonacci numbers:
\begin{equation}\label{3.3}
F^{(1)}_n= F_n,\;\;\;n\geq 0.
\end{equation}

\bfseries Definition 3.\enskip\mdseries  For a positive integer $m$, the
sequence is defined by the same recursion
\begin{equation}\label{3.4}
f^{(m)}_n= f^{(m)}_{n-1} + f^{(m)}_{n-m-1},
\end{equation}
but with other initial conditions

$$f^{(m)}_0=
f^{(m)}_1=\ldots=f^{(m)}_{\lfloor\frac{m-1}{2}\rfloor}=1,\;\;$$
\begin{equation}\label{3.5}
f^{(m)}_{\lfloor\frac{m+1}{2}\rfloor}=f^{(m)}_{\lfloor\frac{m+1}{2}\rfloor+1}=\ldots=
f^{(m)}_m=2,
\end{equation}
we call $m$-Fibonacci numbers of type 2.

For $m=1$ we also obtain the very Fibonacci numbers, but, comparing with
(\ref{3.3})), we have
\begin{equation}\label{3.6}
f_n^{(1)}=F_{n+1},\;\;\;n\geq 0.
\end{equation}
Now we prove several lemmas.

\begin{lemma}\label{L3}
The following  formulas hold
\begin{equation}\label{3.7}
F^{(m)}_n=\sum\limits_{k\geq 0}\begin{pmatrix} n-mk\\
k\end{pmatrix}, \;\;n\geq 0,\;\;m\geq 1,
\end{equation}

\begin{equation}\label{3.8}
f^{(m)}_n=\begin{cases}\sum\limits_{k\geq
0}\begin{pmatrix}n+\frac{m+1}{2}-mk\\k\end{pmatrix}, &\text{if $m$
is odd},\\\sum\limits_{k\geq 0}\begin{pmatrix}n+\frac m
2-(m-1)k\\k\end{pmatrix}, &\text{if $m$ is even}.\end{cases}
\end{equation}
In particular, for odd $m$ we have
\begin{equation}\label{3.9}
f^{(m)}_n= F^{(m)}_{n+\frac{m+1}{2}}.
\end{equation}
\end{lemma}

\bfseries Proof\mdseries\;\; of the (\ref{3.7})-(\ref{3.8}) is
over by the same scheme. Therefore, we prove (\ref{3.7}) only.
Denote the right part of (\ref{3.7}) via $\Phi^{(m)}_n.$  We have
$$
\Phi^{(m)}_{n-1}+\Phi^{(m)}_{n-m-1}=$$
 $$\sum\limits_{k\geq 0}\begin{pmatrix} n-mk-1\\
k\end{pmatrix}+\sum\limits_{k\geq 0}\begin{pmatrix} n-m(k+1)-1\\
k\end{pmatrix}=$$

$$
\sum\limits_{k\geq 0}\begin{pmatrix} n-mk-1\\
k\end{pmatrix}+\sum\limits_{k\geq 1}\begin{pmatrix} n-mk-1\\
k-1\end{pmatrix}=$$
 $$\sum\limits_{k\geq 0}\begin{pmatrix} n-mk\\
k\end{pmatrix}=\Phi^{(m)}_n.
$$

Besides, for $k\geq 1,\;\;m\leq n,$ we have $n-mk\leq m-mk\leq 0.$
\newpage
Hence,
$$
\sum\limits_{k\geq 0}\begin{pmatrix} n-mk\\
k\end{pmatrix}=\begin{pmatrix} n-mk\\
k\end{pmatrix}=1, \enskip if \enskip k=0.
$$

Thus, for numbers $\Phi^{(m)}_n$ and $F^{(m)}_n,$  formulas (\ref{3.1})-(\ref{3.2})
are valid and, consequently, $\Phi^{(m)}_n=F^{(m)}_n,\;\;n\geq
0.$\newline
$\blacksquare$

\begin{lemma}\label{L4}
Let $h_k\equiv k\pmod{2},\enskip h_k=0 \enskip or \enskip 1,$ and $\gamma_m\equiv
m-1\pmod{2},\;\;\gamma_m=0 \enskip or \enskip1.$
\end{lemma}
Then
\begin{equation}\label{3.10}
\sum\limits_{k\geq
1}\begin{pmatrix}\left\lfloor\frac{n-mk-h_k}{2}\right\rfloor\\
\left\lfloor\frac k 2\right\rfloor\end{pmatrix}=
f^{(m)}_{\left\lfloor\frac{n-\gamma_m}{2}\right\rfloor}-1.
\end{equation}

\bfseries Proof\mdseries\;\;\;is the same for all
assumptions regarding the parity of $m$ and $n$. Therefore, we consider only case when $m$ and $n$ are even. Then $\gamma=1$
and
$$
\left\lfloor\frac{n-mk-h_k}{2}\right\rfloor=\frac{n-mk}{2}-h_k.
$$

Therefore, it is enough to prove that
\begin{equation}\label{3.11}
\sum\limits_{k\geq
0}\begin{pmatrix}\frac 1 2(n-mk)-h_k\\
\left\lfloor\frac k 2\right\rfloor\end{pmatrix}= f^{(m)}_{\frac n
2-1}.
\end{equation}

Denote the left part of (\ref{3.11}) via $G^{(m)}_{\frac n 2}$. First
of all, let us show that $G^{(m)}_{\frac n 2}$ satisfies to
(\ref{3.4}). We have
$$
G^{(m)}_{\frac n 2-1}+G^{(m)}_{\frac n 2-m-1}=$$ $$\sum\limits_{k\geq 0}
\begin{pmatrix}\frac n 2 -\frac{mk+2h_k}{2}-1\\\lfloor\frac k 2\rfloor
\end{pmatrix}+\sum\limits_{k\geq 0}
\begin{pmatrix}\frac n 2 -\frac{m(k+2)+2h_{k+2}}{2}-1\\\lfloor\frac k 2\rfloor
\end{pmatrix}=
$$
$$
=\sum\limits_{k\geq 0}
\begin{pmatrix}\frac n 2 -\frac{mk+2h_k}{2}-1\\\lfloor\frac k 2\rfloor
\end{pmatrix}+\sum\limits_{k\geq 2}
\begin{pmatrix}\frac n 2 -\frac{mt+2h_t}{2}-1\\\lfloor\frac t
2\rfloor-1\end{pmatrix}=
$$
$$
=2+\sum\limits_{k\geq 2}
\begin{pmatrix}\frac n 2 -\frac{mk+2h_k}{2}\\\lfloor\frac k 2\rfloor
\end{pmatrix}=\sum\limits_{k\geq 0}
\begin{pmatrix}\frac n 2 -\frac{mk+2h_k}{2}\\\lfloor\frac k 2\rfloor
\end{pmatrix}=G^{(m)}_{\frac n 2}.
$$

It is left to verify the coincidence of the initial conditions (\ref{3.5}) for
$f^{(m)}_{i-1}$ and for
\begin{equation}\label{3.12}
G^{(m)}_i =\sum\limits_{k\geq
0}\begin{pmatrix}i-\frac{mk}{2}-h_k\\\lfloor\frac k
2\rfloor\end{pmatrix}\;\;i=1,2,\ldots, m.
\end{equation}
Notice that, the summands in (\ref{3.12}) equal to $0$ for those and
only those $k$ for which the following inequality is satisfied:
$$
i-\frac{mk}{2}\geq\begin{cases}\frac k 2, \;\; if\;\;k \;\; is \;\;even,\\
\frac{k+1}{2}, \;\; if\;\;k \;\; is\;\; odd.\end{cases}
$$\newpage
or, the same, for
$$
k\geq\begin{cases}\frac{2i}{m+1}, \;\; if\;\;k \;\; is \;\;even,\\
\frac{2i-1}{m+1}, \;\; if\;\;k \;\; is\;\; odd.\end{cases}
$$

Thus, for $i=1,2,\ldots,\frac m 2,$  we have $k=0.$ Thus, by (\ref{3.12}),
$$G^{(m)}_i=\begin{pmatrix} i \\ 0 \end{pmatrix}=1$$

and, for $i=\frac m 2 +1,\;\;\frac m 2 +2,\ldots, m$ we have $k=0$ or
$k=1$.

Therefore, by (\ref{3.12}), $$G^{(m)}_i=\begin{pmatrix} i\\
0\end{pmatrix}+\begin{pmatrix} i-\frac m 2 -1\\
0\end{pmatrix}=2$$

and the lemma follows. \newline
$\blacksquare$

\begin{remark}
\slshape Let show that there does not exist an extension of a
definition of $\left\{f^{(m)}_n\right\}$ for $m=0$ such that the equality
$(\ref{3.10})$ holds.\upshape

Indeed, for $m=0,$ the left part of (\ref{3.10}) equals to
$$\sum\limits_{k=2,4,\ldots.}\begin{pmatrix}\lfloor\frac n
2\rfloor\\\frac k
2\end{pmatrix}+\sum\limits_{k=1,3,\ldots.}\begin{pmatrix}\lfloor\frac
{n-1}{2}\rfloor\\\frac {k-1}{2}\end{pmatrix},$$ while the right part
equals to "$f^{(n)}_{\lfloor\frac{n-1}{2}\rfloor}-1".$ This means
that there must exist a function, say, $g(x)$ such that
$$
\sum\limits_{k=2,4,\ldots}\begin{pmatrix}\lfloor\frac n 2\rfloor\\
\frac k
2\end{pmatrix}=g\left(\left\lfloor\frac{n-1}{2}\right\rfloor\right),
$$
or
$$
g\left(\left\lfloor\frac{n-1}{2}\right\rfloor\right)=
\begin{cases} 2^{\frac n 2}-1,\;\;if\;\;n\;\;is\;\;even,\\
2^{\frac{n-1}{2}}-1,\;\;if\;\;n\;\;is\;\;odd.\end{cases}
$$
For $n=3$ and $n=4,$ we simultaneously have  $g(1)=1$ and $g(1)=3.$\newline
 $\blacksquare$
\end{remark}

\begin{lemma}\label{L5}
The following formula holds

$$\sum\limits^{\lfloor\frac{n}{m+1}\rfloor}_{i=1}\frac 1
k\sum\limits_{d|(n,k)}\varphi(d)\begin{pmatrix}\frac{n-mk}{d}-1\\
\frac k d -1\end{pmatrix}=$$
\begin{equation}\label{3.13}
\frac 1 n
\sum\limits_{\substack{d|n\\d\leq\lfloor\frac{n}{m+1}\rfloor}}
\varphi(d)\left((m+1)F^{(m)}_{\frac n d}-mF^{(m)}_{\frac n
d}-1\right).
\end{equation}
\end{lemma}

\bfseries Proof.\mdseries\;\;\; Put $k=dd_1$. Then we have
$$
\sum\limits^{\lfloor\frac{n}{m+1}\rfloor}_{k=1}\frac 1
k\sum\limits_{d|(n,k)}\varphi(d)\begin{pmatrix}\frac{n-mk}{d}-1\\
\frac k d -1\end{pmatrix}=$$\newpage $$ \sum\limits_{\substack{
dd_1\leq\lfloor\frac{n}{m+1}\rfloor\\d|n}}\frac{1}{dd_1}\;\varphi(d)
\begin{pmatrix}\frac n d -md_1-1\\d_1-1\end{pmatrix}=
$$
\begin{equation}\label{3.14}
=\sum\limits_{\substack{d\leq
\lfloor\frac{n}{m+1}\rfloor\\d|n}}\;\frac{\varphi(d)}{d}\sum\limits^
{\lfloor\frac{n}{(m+1)d}\rfloor}_{d_1=1}\;\frac{1}{d_1}\begin{pmatrix}
\frac n d -md_1-1\\d_1-1\end{pmatrix}.
\end{equation}
Put $h=\frac n d$. Notice that $d\leq\frac{n}{m+1}$, therefore
$t\geq m+1$. \newline
Furthermore,
$$
\sum\limits_{d_1\geq
1}\frac{h}{d_1}\begin{pmatrix}h-md_1-1\\d_1-1\end{pmatrix}=$$ $$\sum\limits_{d_1\geq
1}\frac{h-md_1}{d_1}\begin{pmatrix}h-md_1-1\\d_1-1\end{pmatrix}+$$

$$+m\sum\limits_{d_1\geq
1}\begin{pmatrix}h-md_1-1\\d_1-1\end{pmatrix}=$$ $$ \sum\limits_{d_1\geq
1}\begin{pmatrix}h-md_1\\d_1\end{pmatrix}+ m \sum\limits_{k\geq
0}\begin{pmatrix}(h-m-1)-mk\\k\end{pmatrix}.
$$
Thus, by Lemma \ref{L3}
$$
\sum\limits_{d_1\geq 1}\frac{1}{d_1}\begin{pmatrix}\frac n d
-md_1-1\\d_1-1\end{pmatrix}=\frac d n\left(F^{(m)}_{\frac n
d}+mF^{(m)}_{\frac n d-m-1}-1\right).
$$
Substituting this to (\ref{3.14}) and taking in to account that
$$
F_{\frac n d-m-1}=F_{\frac n d}-F_{\frac n d -1},
$$
we obtain (\ref{3.13}). \newline
$\blacksquare$.

\section{Variations 5-8}

Let now consider two-color bracelets with $k$ isolated black beads such that between every two black beads there exist at least $m$ white ones. Let us  consider an aggregate, denoting it $ black^*,$ which
contains a black bead and the following in succession after it in a fixed
direct the $m$ 0`s. This gives a one-to-one
correspondence between the considered set of all two-color bracelets with $k$ isolated black beads such that between every two black beads there exist at least $m$ white ones and the set of all
two-color bracelets of
$n-mk$ beads containing $k$ $black^*$ beads. Therefore, it follows from
Theorem \ref{t1} that there are $R(n-mk, \;k)$ incongruent
configurations. Notice that always $n-mk\geq k$, or $k\leq
\lfloor\frac{n}{m+1}\rfloor$.

Summing $R(n-mk,\;k)$ by $k$ from 1 to $\lfloor\frac{n}{m+1}\rfloor$
we find the number $\alpha^{(m)}_n$\newpage of all two-color incongruent bracelets with isolated black beads such that between every two black beads there exist at least $m$ white ones.

\begin{equation}\label{4.1}
\alpha^{(m)}_n=\frac 1 2
\sum\limits^{\lfloor\frac{n}{m+1}\rfloor}_{k=1} \left(\frac 1 k
\sum\limits_{d|(k,n)}\varphi(d)\begin{pmatrix}\frac{n-mk}{d}-1\\\frac
k d
-1\end{pmatrix}+\begin{pmatrix}\lfloor\frac{n-mk-h_k}{2}\rfloor\\\lfloor\frac
k 2\rfloor\end{pmatrix}\right),
\end{equation}
where $h_k\equiv k(\mod{2}),\;\;h_k=0\; or \;1,\;\; m\geq 1$.

By (\ref{4.1}) and Lemmas \ref{L4}, \ref{L5} we have

$$\alpha^{(m)}_n=\frac 1
2f^{(m)}_{\lfloor\frac{n-\gamma_m}{2}\rfloor}-\frac 1 2+$$
\begin{equation}\label{4.2}
\frac {1}{2n}
\sum\limits_{d|n,\enskip d\leq\lfloor\frac{n}{m+1}\rfloor}\varphi(d)\left((m+1)F^{(m)}_{\frac
n d}-mF^{(m)}_{\frac n d}-1\right),
\end{equation}
where $\gamma_m\equiv m-1(\mod{2}),\;\;\gamma_m=0\;\;or\;\;1$.\newline
 Now (as in proof of (\ref{2.2})) we notice that the binomial coefficient
$\begin{pmatrix}\lfloor\frac{n-mk-h_k}{2}\rfloor\\\lfloor\frac k
2\rfloor\end{pmatrix}$ gives the number of those
two-color  bracelets with $k$ black beads of the considered type
that are symmetric respectively any diameter. Thus,
by Lemma \ref{L4}, the number $\beta^{(m)}_n$ of all these
bracelets equals to
\begin{equation}\label{4.3}
\beta^{(m)}_n=f^{(m)}_{\lfloor\frac{n-\gamma_m}{2}\rfloor}-1,\;\;\gamma_m\equiv
m-1 \pmod{2},\;\;\gamma_m=0\; or\; 1.
\end{equation}

Now,  using the scheme of the proof of Theorem \ref{t2}, we obtain
the following results.
\begin{theorem}\label{t3}
The following formulas solve Variations $\ref{V5}-\ref{V8}:$
\begin{equation}\label{4.4}
N_{n,m}=-\sum\limits_{\delta |n,\enskip\delta\geq 2}
\mu(\delta){\alpha}^{(m)}_{\frac n \delta},
\end{equation}
where
$\alpha^{(m)}_n$ is defined by $(\ref{4.2});$
\begin{equation}\label{4.5}
 S_{n,m}=-1-\sum\limits_{\delta |n,\enskip\delta\geq 2}\mu(d)\beta^{(m)}_{\frac n d},
\end{equation}
where $\beta^{(m)}_n$ is defined by $(\ref{4.3});$
\begin{equation}\label{4.6}
N^{(k)}_{n,m}=-\sum\limits_{d|(n,k),\enskip d\geq
2}\mu(d)R\left(\frac{n-mk}{d},\;\frac k d\right);
\end{equation}
\begin{equation}\label{4.7}
S^{(k)}_{n,m}=-\sum\limits_{d|(n,k),\enskip d\geq
2}\mu(d)\begin{pmatrix}\left\lfloor\frac 1
2\left(\frac{n-mk}{d}-h_{\frac k
d}\right)\right\rfloor\\\left\lfloor\frac{k}{2d}\right\rfloor\end{pmatrix}.
\end{equation}
\end{theorem}
\;\;\;\;\;

\section{An example}

In the case of $m=1,$ by (\ref{3.3}) and (\ref{3.6}), we have
$$
F^{(1)}_n=F_n,\;\;f^{(1)}_n=F_{n+1},
$$

where
$$
F_0=1,\;F_1=1,\;F_2=2,\;F_3=3,\;F_4=5,\;F_5=8,\;F_6=13,\ldots
$$
are the Fibonacci numbers.

Thus, by (\ref{4.2})-(\ref{4.4}),

$$
{\alpha}^{(1)}_n=\frac{1}{2n}\sum\limits_{d|n,d\leq
\lfloor\frac n 2\rfloor}\varphi(d)\left(2F_{\frac n d}-F_{\frac n
d}-1\right)+\frac 1 2 F_{\lfloor\frac n 2\rfloor+1}-\frac{1}{2},
$$
$$
\beta^{(1)}_n=F_{\lfloor\frac n 2\rfloor+1}-1.
$$

Therefore, by Theorem \ref{t3}, we find
$$
N_{24,1}=
{\alpha}^{(1)}_{12}+{\alpha}^{(1)}_{8}-{\alpha}^{(1)}_{4}=
25+7-2=30,
$$
$$
S_{24,1}=-1+F_{\frac{24}{2\cdot 2}+1}+F_{\frac{24}{2\cdot
3}+1}-F_{\frac{24}{2\cdot 6}+1}=-1+21+8-3=25.
$$
 Thus among $30$ incongruent two-color $SR$-bracelets of 24 beads
 with isolated black beads only five ones do not have a
 diameter of symmetry.

 Further, by Theorem \ref{t3}, we find

$$
N^{(6)}_{24,1}=-\sum\limits_{d=2,3,6}\mu(d)R\left(\frac{18}{d},\;\frac
6 d\right)=9,
$$

$$
S^{(6)}_{24,1}=-\sum\limits_{d=2,3,6}\mu(d)\begin{pmatrix}\left\lfloor\frac
1 2\left(\frac{18}{d}-h_{\frac 6
d}\right)\right\rfloor\\\lfloor\frac 3 d\rfloor\end{pmatrix}=6,
$$

$$
N^{(8)}_{24,1}=-\sum\limits_{d=2,4,8}\mu(d)R\left(\frac{16}{d},\;\frac
8 d\right)=8,
$$

$$
S^{(8)}_{24,1}=-\sum\limits_{d=2,4,8}\mu(d)\begin{pmatrix}\left\lfloor\frac
1 2\left(\frac{16}{d}-h_{\frac 8
d}\right)\right\rfloor\\\lfloor\frac 4 d\rfloor\end{pmatrix}=6.
$$

Since
$$(N^{(6)}_{24,1}-S^{(6)}_{24,1})+(N^{(8)}_{24,1}-S^{(8)}_{24,1})=5, $$
then only in cases $k=6$ and $k=8$ there are correspondingly $3$
and $2$ two-color $SR$-bracelets of 24 beads
 with isolated black beads which do not have a
diameter of symmetry; \enskip in all other cases they do have a diameter of
symmetry.

In Table 1 we show all $30$ incongruent two-color $SR$-bracelets of 24 beads with isolated black beads. 0's denote white beads and 1's denote black ones.
\newpage
\begin{tab}\label{tab1}
\begin{small}
$$
\begin{matrix}
&ordinal \;number&two-color \enskip $SR$-bracelets&number\enskip of \enskip black\enskip beads\\[8pt]
&1&1\;0\;0\;0\;0\;0\;0\;0\;0\;0\;0\;0\;1\;0\;0\;0\;0\;0\;0\;0\;0\;0\;0\;0&2\\
&2&1\;0\;0\;0\;0\;0\;0\;0\;1\;0\;0\;0\;0\;0\;0\;0\;1\;0\;0\;0\;0\;0\;0\;0&3\\
&3&0\;0\;0\;1\;0\;0\;0\;0\;0\;1\;0\;0\;0\;0\;0\;1\;0\;0\;0\;0\;0\;1\;0\;0&4\\
&4&0\;0\;0\;0\;0\;1\;0\;1\;0\;0\;0\;0\;0\;0\;0\;0\;0\;1\;0\;1\;0\;0\;0\;0&4\\
&5&0\;1\;0\;0\;0\;0\;0\;0\;0\;0\;1\;0\;0\;1\;0\;0\;0\;0\;0\;0\;0\;0\;1\;0&4\\
&6&0\;0\;1\;0\;0\;0\;0\;0\;0\;0\;1\;0\;0\;0\;1\;0\;0\;0\;0\;0\;0\;0\;1\;0&4\\
&7&0\;0\;0\;1\;0\;0\;0\;0\;1\;0\;0\;0\;0\;0\;0\;1\;0\;0\;0\;0\;1\;0\;0\;0&4\\
&8&0\;0\;1\;0\;0\;0\;1\;0\;0\;0\;1\;0\;0\;0\;1\;0\;0\;0\;1\;0\;0\;0\;1\;0&6\\
&9&0\;0\;0\;1\;0\;1\;0\;0\;0\;0\;0\;1\;0\;1\;0\;0\;0\;0\;0\;1\;0\;1\;0\;0&6\\
&10&0\;1\;0\;0\;0\;0\;1\;0\;0\;1\;0\;0\;0\;0\;1\;0\;0\;1\;0\;0\;0\;0\;1\;0&6\\
&11&0\;0\;0\;0\;1\;0\;1\;0\;1\;0\;0\;0\;0\;0\;0\;0\;1\;0\;1\;0\;1\;0\;0\;0&6\\
&12&1\;0\;0\;0\;0\;0\;0\;1\;0\;0\;1\;0\;1\;0\;0\;0\;0\;0\;0\;1\;0\;0\;1\;0&6\\
&13&1\;0\;0\;0\;0\;0\;1\;0\;0\;0\;1\;0\;1\;0\;0\;0\;0\;0\;1\;0\;0\;0\;1\;0&6\\
&14&0\;0\;0\;1\;0\;0\;1\;0\;0\;1\;0\;0\;0\;0\;0\;1\;0\;0\;1\;0\;0\;1\;0\;0&6\\
&15&1\;0\;0\;0\;0\;1\;0\;1\;0\;0\;0\;0\;1\;0\;0\;0\;0\;1\;0\;1\;0\;0\;0\;0&6\\
&16&1\;0\;0\;0\;0\;1\;0\;0\;0\;1\;0\;0\;1\;0\;0\;0\;0\;1\;0\;0\;0\;1\;0\;0&6\\
&17&1\;0\;0\;1\;0\;0\;1\;0\;0\;1\;0\;0\;1\;0\;0\;1\;0\;0\;1\;0\;0\;1\;0\;0&8\\
&18&0\;0\;0\;1\;0\;1\;0\;1\;0\;1\;0\;0\;0\;0\;0\;1\;0\;1\;0\;1\;0\;1\;0\;0&8\\
&19&1\;0\;0\;0\;0\;1\;0\;0\;1\;0\;1\;0\;1\;0\;0\;0\;0\;1\;0\;0\;1\;0\;1\;0&8\\
&20&0\;0\;1\;0\;1\;0\;0\;1\;0\;1\;0\;0\;0\;0\;1\;0\;1\;0\;0\;1\;0\;1\;0\;0&8\\
&21&1\;0\;0\;0\;1\;0\;1\;0\;1\;0\;0\;0\;1\;0\;0\;0\;1\;0\;1\;0\;1\;0\;0\;0&8\\
&22&1\;0\;0\;0\;1\;0\;0\;1\;0\;0\;1\;0\;1\;0\;0\;0\;1\;0\;0\;1\;0\;0\;1\;0&8\\
&23&0\;0\;1\;0\;0\;1\;0\;1\;0\;0\;1\;0\;0\;0\;1\;0\;0\;1\;0\;1\;0\;0\;1\;0&8\\
&24&0\;0\;1\;0\;1\;0\;0\;0\;1\;0\;1\;0\;0\;0\;1\;0\;1\;0\;0\;0\;1\;0\;1\;0&8\\
&25&1\;0\;1\;0\;0\;0\;1\;0\;1\;0\;1\;0\;0\;0\;1\;0\;1\;0\;1\;0\;0\;0\;1\;0&9\\
&26&1\;0\;0\;1\;0\;1\;0\;0\;1\;0\;0\;1\;0\;1\;0\;0\;1\;0\;0\;1\;0\;1\;0\;0&9\\
&27&1\;0\;1\;0\;1\;0\;0\;0\;1\;0\;1\;0\;1\;0\;1\;0\;1\;0\;0\;0\;1\;0\;1\;0&10\\
&28&1\;0\;0\;1\;0\;1\;0\;1\;0\;1\;0\;0\;1\;0\;0\;1\;0\;1\;0\;1\;0\;1\;0\;0&10\\
&29&1\;0\;1\;0\;0\;1\;0\;1\;0\;0\;1\;0\;1\;0\;1\;0\;0\;1\;0\;1\;0\;0\;1\;0&10\\
&30&1\;0\;1\;0\;1\;0\;1\;0\;1\;0\;1\;0\;1\;0\;1\;0\;1\;0\;1\;0\;1\;0\;1\;0&12
\end{matrix}
$$
\end{small}
\end{tab}

Note that in Table 1 only bracelets with the ordinal
numbers 12,13,16,19 and 22 have no a diameter of symmetry; it is interesting that the
bracelets with the ordinal numbers 5,7,10,20 do not
have a diameter of symmetry that connects any two beads. More exactly, a diameter of symmetry of these
bracelets has the following endpoints: the midpoint between
positions 24,1 and the midpoint between positions 12,13. A diameter
of symmetry of other bracelets connects beads 1 and 13.

\section{Other variations}
We start with the two open questions arising as
observations of Table \ref{tab1}.

\begin{variation} \label{V9} To enumerate incongruent two-color $SR$-bracelets
with a diameter of symmetry connecting any two beads.\upshape
\end{variation}
Notice, furthermore, that among 21 such bracelets in Table \ref{tab1} only 3 ones
\newpage
have a diameter of symmetry connecting beads of different colors (see
bracelets with ordinal numbers 2, 25, 26 ). The following enumerating problem arises naturally.
\begin{variation} \label{V10} To enumerate incongruent two-color $SR$-bracelets
with a diameter of symmetry connecting any two beads of different colors.
\end{variation}

In conclusion, we give several cases of enumeration of the incongruent
$t$-color bracelets $(t\geq3)$ containing beads of every of $t$ colors.
\begin{variation} \label{V11}
To find the number
$N^{(0,1,2)}_{n,1}$ incongruent three-color $($with colors $\{0,1,2\})$ $SR$-bracelets
with isolated beads of colors $1$ and $2.$
\end{variation}

\begin{theorem}\label{t4}

\begin{equation}\label{6.1}
N^{(0,1,2)}_{n,1}=\sum\limits^{\lfloor\frac n
2\rfloor}_{\substack{k\geq2\\(k,n)>1}}N^{(k)}_{n,1}N_{(k,n)}.
\end{equation}
\end{theorem}
\bfseries Proof.\mdseries  \enskip We generate the required bracelets from $N^{(k)}_{n,1}$ two-color $SR$-bracelets with colors 0 and 1, $k,\enskip 2\leq k\leq \lfloor(\frac{n}{2})\rfloor.$ Removing beads of color 0, we obtain $k$ places for generating two-color subbracelet of colors 1 and 2. Since we consider $SR$-bracelets, then, in view of  subgroup of subgroup of a cyclic grope is a subgroup, we  should have $(l,k)|(k,n)|n.$ Thus
\begin{equation}\label{6.2}
N^{(0,1,2)}_{n,1}=\sum\limits^{\lfloor\frac n
2\rfloor}_{\substack{k\geq2\\(k,n)>1}}N^{(k)}_{n,1} \sum\limits^{\lfloor\frac {(k,n)}
2\rfloor}_{\substack{l\geq2\\1<(l,k)|(k,n)}}N^{(l)}_{(k,n)}.
\end{equation}
Since
$$ \sum\limits^{\lfloor\frac n
2\rfloor}_{\substack{l\geq2\\(l,n)>1}}N^{(l)}_{n}=N_n$$
(and always $(l,n)|n),$ then from (6.2) we obtain (\ref{6.1}).
\newline
$\blacksquare$\newline

Thus the required enumeration we get from  formulas (\ref{6.1}), (\ref{2.1}), (\ref{4.4}) and (\ref{4.6}).

Let us calculate, e.g., $N^{(0,1,2)}_{12,1}.$ By (\ref{6.1}), we have
$$N^{(0,1,2)}_{12,1}=N_{12,1}^{(2)}N_2+N_{12,1}^{(3)}N_3+N_{12,1}^{(4)}N_4+N_{12,1}^{(6)}N_6.  $$
Since, for prime $p,\enskip N_p=0,$ then
\begin{equation}\label{6.3}
N^{(0,1,2)}_{12,1}=N_{12,1}^{(4)}N_4+N_{12,1}^{(6)}N_6.
\end{equation}
According to (\ref{2.5})-(\ref{2.6}), we find $\alpha_i=i+1, \enskip i=1,2,3.$ Thus, by (\ref{2.1}), \newpage we have
$$N_4=-2+\alpha_2=1,\enskip N_6=-2+\alpha_3+\alpha_2-\alpha_1=3.$$

Furthermore, by (\ref{4.6}) and (\ref{1.1}), we have
$$N_{12,1}^{(4)}=R(4,2)=2  $$
and
$$N_{12,1}^{(6)}=R(3,3)+R(2,2)-R(1,1).  $$
Since, by (\ref{1.1}),
 $$R(n,n)=\frac{1}{2}(1+\frac{1}{n}\sum_{d|n}\varphi(d))=1,$$
 then $N_{12,1}^{(6)}=1$ and, by (\ref{6.3}), finally we find
$$ N^{(0,1,2)}_{12,1}=2\cdot1+1\cdot3=5.  $$

In Table \ref{tab2} we demonstrate all 5 these three-color bracelets of 12 beads.

\begin{tab}\label{tab2}

\begin{small}
$$
\begin{matrix}
&ordinal \;number&3-color\enskip bracelets\\[8pt]

&1&1\;0\;0\;2\;0\;0\;1\;0\;0\;2\;0\;0\\
&2&1\;0\;1\;0\;2\;0\;1\;0\;1\;0\;2\;0\\
&3&2\;0\;2\;0\;1\;0\;2\;0\;2\;0\;1\;0\\
&4&1\;0\;2\;0\;1\;0\;2\;0\;1\;0\;2\;0\\
&5&1\;0\;2\;0\;0\;0\;1\;0\;2\;0\;0\;0
\end{matrix}
$$
\end{small}
\end{tab}
\begin{variation} \label{V12}
To find the number
$S^{(0,1,2)}_{n,1}$ incongruent three-color $($with colors $\{0,1,2\})$ $SR$-bracelets with a diameter of symmetry
having isolated beads of colors $1$ and $2.$
\end{variation}

A solution is given by the following similar to (\ref{6.1})
formula.
\begin{theorem}\label{t5}
\begin{equation}\label{6.4}
S^{(0,1,2)}_{n,1}=\sum\limits^{\lfloor\frac n
2\rfloor}_{\substack{k=2\\(k,n)>1}}S^{(k)}_{n,1} S_{(k,n)}.
\end{equation}
\end{theorem}
Proof of Theorem \ref{t5} is based on quite similar arguments that proof of Theorem \ref{t4}.
Thus the required enumeration we obtain from  formulas (\ref{6.2}), (\ref{2.2}), (\ref{4.5}) and (\ref{4.7}).
So, e.g., we find that $S^{(0,1,2)}_{12,1}=14.$ Thus there is only (up to the congruence) three-color $SR$-bracelet in Table \ref{tab2} that has not a diameter of symmetry. It is easy to see that it is the last bracelet in this table.\newline
\indent Finally note that the idea of construction formulas (\ref{6.1}), (\ref{6.4}) is easily generalized. E. g., if  $N^{(0,1,...,t)}_{n}\enskip (S^{(0,1,...,t)}_{n})$ denotes the number of incongruent $(t+1)$-color $SR$-bracelets
(having a diameter of symmetry),
\newpage
 $N^{(0,1,...,t)}_{n,1}\enskip (S^{(0,1,...,t)}_{n,1})$ denotes the number of incongruent $(t+1)$-color $SR$-bracelets
with isolated beads of colors $1,...,t$ (having a diameter of symmetry), etc., then we have the following recursion formulas.

\begin{theorem}\label{t6}
\begin{equation}\label{6.5}
N^{(0,1,...,t)}_{n}=\sum\limits^{\lfloor\frac n
2\rfloor}_{\substack{k=2\\(k,n)>1}}N^{(k)}_{n} N^{(0,1,...,t-1)}_{(k,n)}.
\end{equation}
\begin{equation}\label{6.6}
S^{(0,1,...,t)}_{n}=\sum\limits^{\lfloor\frac n
2\rfloor}_{\substack{k=2\\(k,n)>1}}S^{(k)}_{n} S^{(0,1,...,t-1)}_{(k,n)}.
\end{equation}

\begin{equation}\label{6.7}
N^{(0,1,...,t)}_{n,1}=\sum\limits^{\lfloor\frac n
2\rfloor}_{\substack{k=2\\(k,n)>1}}N^{(k)}_{n,1} N^{(0,1,...,t-1)}_{(k,n)}.
\end{equation}
\begin{equation}\label{6.8}
S^{(0,1,...,t)}_{n,1}=\sum\limits^{\lfloor\frac n
2\rfloor}_{\substack{k=2\\(k,n)>1}}S^{(k)}_{n,1} S^{(0,1,...,t-1)}_{(k,n)}.
\end{equation}
etc.
\end{theorem}


\begin{thebibliography}{6}
\bibitem{1} 1. H.\enskip Gupta,\enskip Enumeration of incongruent cyclic k-gons,
\slshape \enskip Indian J.\enskip Pure Appl. Math.,\upshape \enskip\bfseries 10 \mdseries (1979),\enskip
no.\enskip 8, \enskip964-999.

\bibitem{2} 2. M.\enskip Hall,\enskip Combinatorial Theory,\slshape \enskip Blaisdell,\upshape\enskip
1967.
\bibitem{3} 3. K.\enskip Ireland, M.\enskip Rosen,\enskip A Classical Introduction to
Modern Number Theory.\slshape \enskip Springer-Verlag,\upshape \enskip1982.

\bibitem{4} 4. V.\enskip Shevelev, \enskip Necklaces and convex k-gons,\slshape \enskip Indian
J. Pure Appl. Math., \upshape \enskip\bfseries 35 \mdseries (2004),\enskip no.\enskip5,\enskip
629-638.
\bibitem{5} 5. V.\enskip Shevelev, http://arxiv.org./abs/1104.4051 (2011), Spectrum of permanent's values and its extremal magnitudes in $\Lambda_n^3$ and $\Lambda_n(\alpha,\beta,\gamma).$
\bibitem {6} 6. N.\enskip J.\enskip A.\enskip Sloane,\enskip\slshape The On-Line Encyclopedia of Integer Sequences \upshape \enskip http://oeis.org.
\end{thebibliography}
\end{document}